\documentclass[journal]{IEEEtran}

\usepackage{cite}

\usepackage{amsmath}

\usepackage{mathtools} 
\usepackage{graphicx,color,multirow,colortbl,hhline,array,url,fancyhdr,setspace,enumerate}
\usepackage{subcaption}
\usepackage{enumitem}
\usepackage{amsmath,bm}
\usepackage{mathtools} 
\usepackage{float}

\makeatletter
\def\widebar{\accentset{{\cc@style\underline{\mskip11mu}}}}
\makeatother
\usepackage{amsthm}
\usepackage{amssymb}

\newtheorem*{remark}{Remark}

\newtheorem{innercustomgeneric}{\customgenericname}
\providecommand{\customgenericname}{}
\newcommand{\newcustomtheorem}[2]{%
	\newenvironment{#1}[1]
	{%
		\renewcommand\customgenericname{#2}%
		\renewcommand\theinnercustomgeneric{##1}%
		\innercustomgeneric
	}
	{\endinnercustomgeneric}
}

\newcustomtheorem{customthm}{Theorem}
\newcustomtheorem{customcorollary}{Corollary}
\newcustomtheorem{customcondition}{Condition}
\newcustomtheorem{customlemma}{Lemma}
\newcustomtheorem{customdefinition}{Definition}
\usepackage[ruled,vlined,linesnumbered]{algorithm2e}
\usepackage{booktabs}
\usepackage{tabularx}

\DeclareMathOperator*{\argmin}{arg\,min}
\usepackage{threeparttable}
\begin{document}

	\title{Power System Dispatch with Marginal Degradation Cost of Battery Storage}

	\author{Guannan~He,~\IEEEmembership{Member,~IEEE,}
		Soummya~Kar,~\IEEEmembership{Member,~IEEE,}
		Javad~Mohammadi,~\IEEEmembership{Member,~IEEE,}
		Panayiotis~Moutis,~\IEEEmembership{Senior~Member,~IEEE,}
		and~Jay~F.~Whitacre
		\thanks{This work was partially supported by the US Department of Energy under Grant DEEE0007165.}
		\thanks{Guannan He was with the Department of Engineering and Public Policy and Wilton E. Scott Institute for Energy Innovation, Carnegie Mellon University, Pittsburgh, PA 15213 USA. He is currently with MIT Energy Initiative, Massachusetts Institute of Technology, Cambridge, MA 02139 USA. (e-mail: gnhe@mit.edu)}
		\thanks{Soummya Kar, Javad Mohammadi, and Panayiotis Moutis are with the Department of Electrical and Computer Engineering, Carnegie Mellon University, Pittsburgh, PA 15213 USA. (e-mail: soummyak@andrew.cmu.edu, jmohamma@andrew.cmu.edu, pmoutis@andrew.cmu.edu)}%
		\thanks{Jay F. Whitacre is with Wilton E. Scott Institute for Energy Innovation, the Department of Engineering and Public Policy, and the Department of Materials Science and Engineering in Carnegie Mellon University. (e-mail: whitacre@andrew.cmu.edu)}
	}
	
	
	\maketitle
	
	\begin{abstract}
		Battery storage is essential for the future smart grid. The inevitable cell degradation renders the battery lifetime volatile and highly dependent on battery dispatch, and thus incurs opportunity cost. This paper rigorously derives the marginal degradation cost of battery for power system dispatch. The derived optimal marginal degradation cost is time-variant to reflect the time value of money and the functionality fade of battery and takes the form of a constant value divided by a discount factor plus a term related to battery state of health. In case studies, we demonstrate the evolution of the optimal marginal costs of degradation that corresponds to the optimal long-term dispatch outcome. We also show that the optimal marginal cost of degradation depends on the marginal cost of generation in the grid.
	\end{abstract}
	

	\IEEEpeerreviewmaketitle

	\section{Introduction}
	
	\IEEEPARstart{E}{lectrochemical} energy storage, also known as battery storage, will be a critical component in the future power system with high penetration of renewable energy \cite{chu12}. By providing flexibility and fast responding capability, battery storage can contribute in many ways such as integrating renewable energy \cite{he17}, alleviating congestion \cite{khani16,elliott18}, providing frequency regulation \cite{he16}, enhancing system stability \cite{ortega18}, and so on. 
	
	The dispatch problem of battery storage has been intensively studied. The dispatch of battery with wind farm to improve the overall performance are investigated in \cite{li11,yao12,abdullah15,khatamianfar13}. The microgrid dispatch model with battery are studied in \cite{battistelli17,hu18,ju18,ross18,hug15}. Many studies also implement distributed optimization algorithm in the dispatch model \cite{hug15,cherukuri18,qin18,xing17}. 
	
	Despite the above documented high utility value of battery storage, there exists no explicit or rigorous way to define its marginal operating cost in literature. Unlike traditional thermal generators that have explicit variable generating costs, the real variable operational cost of battery storage is approximately zero as it consumes no fuel. Many studies ignore battery degradation and thus do not impose a marginal degradation cost for battery in dispatch objective function \cite{li11,yao12,abdullah15,ross18,cherukuri18,xing17}, or do not explain or justify why they adopt such a cost/penalty function \cite{khatamianfar13,hug15,qin18}. Most electrochemical energy storage degrades when it is under use, which seems a legitimate cost source. Some studies on battery storage dispatch incorporate a degradation/aging cost that is derived from the capital or replacement cost of battery storage \cite{battistelli17,hu18,ju18,bordin17,shi18-2nQfj,tan19,zhang15-1Q4GT,swierczynski14,farzin16,ortega-vazquez14-j6PKn,xu17,tant13,xu18-factoring,hoke14}. In \cite{he18}, we argue and prove that the capital or replacement cost should not affect the dispatch of battery storage, in general, because the capital/replacement cost is a long-term cost and should not determine the marginal operational cost. However, the proposed framework in \cite{he18} only presents degradation-constrained-form model, whereas a justified cost function for battery storage is still absent. It is imperative to derive such a cost function, if any, because the solution to a dispatch model with wrong objective function will greatly deviate from the optimal dispatch outcome.
	
	In this paper, we have rigorously derived the marginal degradation costs of battery storage that guarantbattery the long-term dispatch optimality of the system. The main technical contributions of this paper include:
	
	1) Proposing a multi-stage intertemporal model (short-term and long-term) for power system dispatch with battery storage. This model considers storage usage and state of health (SOH) as decision and state variable, respectively. The long-term multi-stage dispatch problem is decomposed into short-term dispatch problems with reduced computational complexity.
	
	2) Providing analytical justification for the marginal degradation costs of storage. Specifically, our justification proves that the there exists a specific time series of marginal degradation costs that can be implemented into the short-term power system dispatch problem to get the optimal long-term system dispatch including the optimal battery usage. 
	
	3) Proposing a backward algorithm that simulates and optimizes the battery SOH evolution from the end to the beginning of battery life to solve the intertemporal model, based on the derived optimal evolution of degradation costs.
	
	Mathematically, the derived unit degradation costs consist of two terms: a) a constant value (namely the marginal cost of degradation, MCD) divided by a discount factor to reflect the time value of money; and b) a term that depends on the sensitivity of system objective with respect to storage SOH. The first term represents the opportunity cost of the current battery usage from the decrease of the battery lifetime, while the second term represents the opportunity cost from the functionality degradation of battery. Implementing the optimal degradation costs proposed in this paper could save significantly more system costs than using a degradation cost derived from the capital or replacement cost as previous studies do \cite{hu18,ju18,bordin17,shi18-2nQfj,tan19,zhang15-1Q4GT}. One key implication of this paper as well as \cite{he18} is that the degradation cost of storage is not a given function but an opportunity cost that depends on other agents and equipment in the power system. 
	
	
	This paper is organized as follows. Section II introduces the generic power system dispatch problem. Section III defines degradation-constrained dispatch and long-term optimization problem. Section IV derives the optimal degradation cost of storage based on the optimal storage usage (Theorem \ref{theorem1} and presents the evolution of the degradation cost (Corollary \ref{corollary1}) as well as a backward algorithm. Section V is a case study of microgrid economic dispatch, in which the evolution of the optimal unit degradation costs is presented, and the results of the proposed and existing methods are compared. Section VI draws the conclusion.
	
	\section{Power System Dispatch Problem}
	
	\subsection{Generic Power System Dispatch Problem}
	\label{sec:studytwo:power_system_dispatch}
	The steady-state dispatch problem of power system is usually formulated as an optimization problem to find the best control of the system with respect to an objective function $ f $ (such as total social welfare maximization, total cost minimization, and so on), subject to a set of operating constraints (such as power balance, line flow limit, voltage limit, and so on) and given parameters (such as load demand, equipment status, cost metrics, and so on), as in equation (\ref{3.1}):
	\begin{equation}
	\label{3.1}
	\begin{split}
	&\min\limits_{\bm{x}} f(\bm{x},\bm{v},\bm{p}) \\
	\mathrm{s.t.} \quad &\bm{g}(\bm{x},\bm{v},\bm{p}) \leq 0
	\\ &\bm{h}(\bm{x},\bm{v},\bm{p}) = 0
	\end{split}
	\end{equation}
	where $ \bm{x} $ are the control variables that may include all quantities controllable for system operator, such as power outputs of generators, demand response, and so on;  $ \bm{v} $ represent the state variables such as phasor angles, voltage magnitudes, and so on; $ \bm{p} $ represent the given parameters such as loads, equipment status and specifications, generating costs, and so on;   represent the inequality constraints, such as line flow limits, voltage limits, generation limits, and so on; $ \bm{h} $ represent the equality constraints, which typically includes power balance constraint.
	
	Depending on application scenarios, the dispatch models of power system take different mathematical forms. For example, the direct-current (DC) optimal power flow (OPF) model ignores the reactive power and assumes that the voltage magnitudes and angles of all buses are the same to eliminate non-linear terms; the economic dispatch (ED) model, as the simplest one, ignores transmission constraints and only implements steady-state power balance constraints. In transient stability assessment, differential equations that represent the dynamics of the system are also added to the model but are beyond the focus of this paper.
	
	Many elements in $ \bm{p} $, such as load and renewable generation, are uncertain. To take advantage of the more accurate forecasting information that is available when the real-time operation is approaching, the dispatch decision is usually made day-ahead, intra-day, intra-hour or rollingly across time. Because battery storage has energy constraints that are time-coupling, the dispatch model with battery storage needs to include multiple time steps. We denote the time index by $ t $, the time horizon by $ \Delta t $, and the optimal value of the objective function by $ F_{t} $. For conciseness, we hide the state variables and parameters and merge the inequality and equality constraints (the latter can be represented by the former), and then the dispatch model for time $ t $ becomes equation (\ref{3.2}):
	\begin{equation}
	\label{3.2}
	\begin{split}
	F_{t} = &\min\limits_{\bm{x}_{t}} f_{t}(\bm{x}_{t}) \\
	\mathrm{s.t.} \quad &\bm{g}(\bm{x}_{t}) \leq 0
	\end{split}
	\end{equation}

	\section{Intertemporal Power System Dispatch with Battery Storage}
	In this section, the definitions of the degradation-constrained dispatch with battery, the long-term objective of dispatch with battery, and the long-term dispatch problem with battery are introduced, respectively.
	
	\subsection{Degradation-constrained power system dispatch problem with battery storage}
	\label{sec:studytwo:IOC}
	When battery is employed, degradation is inevitably incurred. As battery degrades, its power and energy capacities and charge/discharge efficiency decrease. To analyze the optimality of dispatch subject to the life-cycle degradation limit of battery, here the degradation-constrained dispatch that implements battery degradation into dispatch problem is proposed.
	
	\begin{customdefinition}{1} \label{definition1} The degradation-constrained dispatch is formulated as follows:
		
		Problem A:
		\begin{align}
		\label{3.3}
		F_{t}(u_{t},H_{t}) = &\min\limits_{\bm{x}_{t}} f_{t}(\bm{x}_{t})
		\\ \label{3.4} \mathrm{s.t.} \quad &\bm{g}_{t}(\bm{x}_{t},H_{t}) \leq \mathbf{0}
		\\  \label{3.5} &d_{t}(\bm{x}_{t}^{\mathrm{S}},H_{t}) \leq u_{t}
		\\  \label{3.6} &\bm{x}_{t}^{\mathrm{S}} \geq \mathbf{0}
		\end{align}
	\end{customdefinition}
	
	In the degradation-constrained dispatch above, namely Problem A, the decision variables are $\bm{x}_{t}$, the control of all controllable quantities in the power system, which include the control of the battery denoted by $\bm{x}_{t}^{\mathrm{S}}$.
	
	In (\ref{3.3}), $ H_{t} $ denotes the state of health (SOH) of battery, which is commonly used as a measure of the storage functionality and cumulative degradation. Typically, SOH is defined as the ratio of the remaining energy capacity to the initial rated energy capacity. For example, if the initial capacity is 1 MW and the current capacity is 900 kWh, the SOH is 90\%. The remaining power capacity and energy efficiency can also be expressed as functions of SOH. The end of battery life is usually defined as when the SOH drops to a criterial percentage, for example, 80\% or 70\%, because the capacities and efficiency will fade much more drastically afterwards.
	
	The degradation of battery in a grid during time $t$ is denoted by $d_{t}\left( \bm{x}_{t}^{\mathrm{S}},H_{t} \right)$ as a function of battery control $ \bm{x}_{t}^{\mathrm{S}} $ (charging and discharging) and the state of health (SOH) of battery $ H_{t} $ at the beginning of time $ t $. The upper limit on the degradation or usage of battery during time $ t $ is denoted by $ u_{t} $ in (\ref{3.5}).
	
	While the constraint $ \bm{x}_{t}^{\mathrm{S}} \geq \mathbf{0} $ in (\ref{3.6}) is not necessary for general battery dispatch problems, this domain limit makes $ d(\bm{x}_{t}^{\mathrm{S}},H_{t}) $ differentiable with respect to $ \bm{x}_{t}^{\mathrm{S}} $. $ d(\bm{x}_{t}^{\mathrm{S}},H_{t}) $ is only piece-wise differentiable if each element of $ \bm{x}_{t}^{\mathrm{S}} $ represent both charging and discharging outputs of battery (for example, the element of $ \bm{x}_{t}^{\mathrm{S}} $ is positive when discharging and negative when charging). The case when $ d(\bm{x}_{t}^{\mathrm{S}},H_{t}) $ is only piece-wise differentiable will also be discussed later, in which this constraint can be relaxed. To satisfy the domain limit, the discharging and charging outputs can be expressed as different non-negative elements in $ \bm{x}_{t}^{\mathrm{S}} $ (for example, $ \bm{x}_{t}^{\mathrm{S}} = \left[ \bm{x}_{t}^{\mathrm{S+}}, \bm{x}_{t}^{\mathrm{S-}}  \right]^{T}  $, where $ \bm{x}_{t}^{\mathrm{S+}} $ is for discharging and $ \bm{x}_{t}^{\mathrm{S-}} $ is for charging). The optimal dispatch solutions cannot have positive discharging and charging outputs at the same time, which only increases degradation without any additional contribution to the system.
	
	The optimal objective of Problem A for time $ t $, $ F_{t} $, is a function of the usage limit for time $ t $ and the SOH of battery. The degradation or usage limit of battery, $ u_{t} $, is exogeneous for the degradation-constrained dispatch. $ F_{t}(u_{t},H_{t}) $ is generally a monotonically decreasing function with respect to $ u_{t} $ and $ H_{t} $, as the optimal solution when $ u_{t}=u^{\ast} $ is always feasible for the problem above when $ u_{t}=u^{\ast} + \Delta u $, where $ \Delta u > 0 $. One sufficient condition for $ F_{t}(u_{t},H_{t}) $ being convex is that $ f_{t}(\bm{x}_{t}) $, $ \bm{g}_{t}(\bm{x}_{t},H_{t}) $, and $ d_{t}(\bm{x}_{t}^{\mathrm{S}},H_{t}) $ are convex with respect to $ \bm{x}_{t} $ \cite{Fiacco1986}. General optimal power flow problems consisting of $ f_{t}(\bm{x}_{t}) $ and $ \bm{g}_{t}(\bm{x}_{t},H_{t}) $ are non-convex but can be linearized (DC OPF) or relaxed (convex relaxation) to convex problems \cite{low14}. $ d_{t}(\bm{x}_{t}^{\mathrm{S}},H_{t}) $ is also convex for common battery chemistries \cite{shi18,xu18}.
	
	\subsection{Long-term power system dispatch problem with battery storage}
	Based on the degradation-constrained dispatch problem, we derive the long-term multi-period power system dispatch probelm with battery storage:
	\begin{customdefinition}{2}\label{definition2} The long-term objective of dispatch with battery is defined as the sum of the present values of all optimal objectives from degradation-constrained dispatches over the battery lifetime, as follows:
		\begin{equation}
		\label{7}
		\begin{split}
		&y = \sum\limits_{t \leq N} \delta_{t}\min\limits_{\bm{x}_{t} \in \mathbf{\Phi}} f_{t}(\bm{x}_{t})
		= \sum\limits_{t \leq N} \delta_{t} F_{t}\left(u_{t},H_{t}\left(H_{0},\sum\limits_{\tau < t}u_{\tau} \right)  \right) 
		\\ &\mathbf{\Phi} =  \left\lbrace \bm{x}_{t}|\bm{g}_{t}(\bm{x}_{t},H_{t}) \leq \mathbf{0},d(\bm{x}_{t}^{\mathrm{S}},H_{t}) \leq u_{t},\bm{x}_{t}^{\mathrm{S}} \geq \mathbf{0}\right\rbrace 
		\end{split}
		\end{equation}
	\end{customdefinition}
	
	\begin{customdefinition}{3}\label{definition3} The long-term optimization problem of dispatch with battery (Problem B) is defined as follows:
		\begin{align}
		Y(H_{0})  = &\min\limits_{u_{t},\bm{\beta}} \sum\limits_{t \leq N} \delta_{t} F_{t}\left(u_{t},H_{t}\left(H_{0},\sum\limits_{\tau < t}u_{\tau} \right)  \right)  	 \label{8} \\ 
		\mathrm{s.t.} \quad &H_{t+1} = H_{t} - \dfrac{u_{t}}{ U(H_{0})}\left(H_{0} - \underline{H} \right) \label{9} \\
		&
		\sum\limits_{t \leq N} \beta_{t} u_{t} \leq U(H_{0}) \label{10}\\ 
		&\beta_{t} \geq \beta_{t+1} \qquad\qquad\qquad\quad \forall t \leq N \label{11}\\ 
		&u_{t} \geq q_{t} \qquad\qquad\qquad\qquad \forall t \leq N \label{12}\\ 
		&u_{t} \leq q_{t} + \beta_{t}U(H_{0}) \qquad\quad \forall t \leq N \label{13}
		\end{align}
	\end{customdefinition}
	
	In (\ref{7}), $ y $ is the long-term objective, and $ \delta_{t} $ is the discount factor for time $ t $. In (\ref{8}), $ Y $ is the optimal value of the long-term objective function $ y $, dependent on the initial battery SOH $ H_{0} $; $ \delta_{t} $ is the discount factor for time $ t $ usually taking the form of $ \dfrac{1}{(1+r)^{\kappa(t)}} $, where $ r $ is the discount rate, and $ \kappa(t) $ is the year number for time $ t $ from the beginning of the battery project. The long-term optimization problem presented above is a multi-stage programming problem. The main decision variable of the problem is the battery usage $u_{t}$ at each stage/time, and the state variable is the SOH.
	
	The SOH evolution is expressed as (\ref{9}), where $ \underline{H} $ represents the SOH when the battery life ends, and $ H_{0} $ is the initial SOH of the optimization horizon under consideration. Based on (\ref{9}), we can further express $ H_{t} $ as a function of the initial SOH $ H_{0} $ and the accumulative battery usage limits before time $ t $, as follows:
	\begin{equation}
	\label{14}
	\begin{split}
	H_{t+1} &= H_{0} - \dfrac{\sum\limits_{\tau \leq t}u_{\tau}}{ U(H_{0})}\left(H_{0} - \underline{H} \right) 
	\end{split}
	\end{equation}
	For conciseness, $ H_{t}\left(H_{0},\sum\limits_{\tau < t}u_{\tau} \right) $ is abbreviated as $ H_{t} $ in the remainder of the paper. The derivative of SOH with respect to any battery usage limit in the past is:
	\begin{equation}
	\label{15}
	\begin{split}
	\dfrac{\partial H_{t} }{\partial u_{\tau}} = -\dfrac{H_{0} - \underline{H} }{ U(H_{0})}\qquad \forall \tau \leq t
	\end{split}
	\end{equation}
	
	In (\ref{10}), the total degradation before the battery life ends is denoted by $ U $, whose unit can be number of cycles or energy throughput (MWh). $ U $ is a function of $ H_{0} $, which is the initial SOH of the optimization horizon under consideration. $ \beta_{t} $ is a binary variable that indicates whether the battery has reached its end of life ($ \beta_{t}=1 $ indicates the battery is within its life at time $ t $, while $ \beta_{t}=0 $ indicates the battery should be retired or replaced before time $ t $). Let us denote the battery lifetime by, $T$, then $ \beta_{t}=1 $ for $ t<=T $ and $ \beta_{t}=0 $ for $ t>T $.
	
	Battery degradation is positive even if the battery stops operation, because of the existence of calendar degradation that is independent of battery cycling or usage and is only dependent on temperature and the state of charge (SOC) of battery. This phenomenon can be expressed as:
	\begin{equation}
	\label{16}
	\begin{split}
	\inf\limits_{\bm{x}_{t}^{\mathrm{S}}}\left[d_{t}\left( \bm{x}_{t}^{\mathrm{S}},H_{t} \right) \right] = d_{t}\left( \mathbf{0},H_{t} \right) = q_{t}
	\end{split}
	\end{equation}	
	where $ q_{t} $ denotes the calendar degradation of the battery during time $ t $. (\ref{16}) directly leads to (\ref{12}). If $ u_{t} = q_{t} $, there is only calendar degradation and no cycling degradation, which implies that the battery is staying idle and $ \bm{x}_{t}^{\mathrm{S}} = \mathbf{0} $. Constraint (\ref{13}) enforces that the battery stops operation after its end of life ($ u_{t} = q_{t} $ when $ \beta_{t}=0 $). For feasibility, $ U \geq q_{t} $. In (\ref{8})-(\ref{13}), $ N $ is a large number that caps the battery lifetime and can be estimated by $ U(H_{0})/q_{t} $.
	
	The degradation modelling above does not make assumptions on the battery chemistry. For various types of battery, they have different $d_{t}\left( \bm{x}_{t}^{\mathrm{S}},H_{t} \right)$, $ q_{t} $, and $ U $. The degradation function $d_{t}\left( \bm{x}_{t}^{\mathrm{S}},H_{t} \right)$  takes an abstract form and can be highly complicated and non-linear to fit empirical degradation data. While temperature may also affect the degradation rate of battery, typical stationary battery systems are equipped with thermal management system to keep the temperature constant and minimize degradation. Even if the temperature is not constant, the effect of temperature can still be modelled in $d_{t}\left( \bm{x}_{t}^{\mathrm{S}},H_{t} \right)$ by updating $d_{t}$ across time for different temperatures, which does not affect any theoretical results that will be presented in the following sections. 
	Since the degradation/usage limit is typically equal to the degradation/usage ($ u_{t} = d(\bm{x}_{t}^{\mathrm{S}},H_{t}) $) in the optimal solutions to Problem B (\ref{8})-(\ref{13}), the battery usage and its upper limit are not differentiated afterwards.
	
	As renewable generation and electricity demand are uncertain in the short-term dispatch (Problem A), $ F_{t}(u_{t},H_{t}) $ can also represent the expected value of the minimum system cost for time $t$. The long-term optimization problem (Problem B) can be framed as a multi-stage stochastic dynamic programming problem, in which context our following results also hold. 
	
%
	
	\section{Marginal Degradation Cost of battery}
	\label{sec:studytwo:degradation_cost}
	In this section, the optimal marginal degradation cost of battery is presented based on the optimal battery usage, and a backward algorithm is proposed to find the optimal value of degradation cost. Theorem \ref{theorem1} indicates that there exists a degradation cost that can be implemented into short-term dispatch problem to get the optimal battery usage as well as system dispatch. The derived battery cost is essentially an opportunity cost from the loss of future benefit opportunity due to degradation and is independent of the initial capital cost of battery. Then the evolution of the optimal unit degradation costs is revealed in Corollary \ref{corollary1}, which leads to a backward algorithm to solve the whole problem.
	
	\subsection{Optimal marginal degradation cost}
	First, two basic conditions \ref{condition1} and \ref{condition2} are presented for Theorem \ref{theorem1}:
	\begin{customcondition}{1}
		\label{condition1}
		$  F_{t}\left(u_{t},H_{t} \right) $ is piece-wise differentiable with respect to $ u_{t} $.
	\end{customcondition}
	\begin{customcondition}{2}
		\label{condition2}
		$ F_{t}\left(u_{t},H_{t} \right) $ is piece-wise convex with respect to $ u_{t} $.
	\end{customcondition}
	
	\begin{customthm}{1}
		\label{theorem1}
		If Conditions \ref{condition1} and \ref{condition2} hold, there exists a non-negative constant $ c^{\ast} $, such that the optimal solution $ u_{t}^{\ast} $ to Problem B (\ref{8})-(\ref{13}) should satisfy either of (\ref{3.14}) and (\ref{3.15}):
		\begin{align}
		u_{t}^{\ast} = q_{t}
		\label{3.14}
		\end{align}	
		\begin{align}
		u_{t}^{\ast} = d(\bm{x}_{t}^{\mathrm{S}\ast},H_{t}^{\ast})
		\label{3.15}
		\end{align}
		where $ \bm{x}_{t}^{\mathrm{S}\ast} $ is the optimal solution to Problem C, as presented below:
		
		Problem C:
		\begin{equation}
		\label{3.16}
		\begin{split}
		J_{t}(c_{t}^{\ast},H_{t}^{\ast}) = &\min\limits_{\bm{x}_{t}} \left[ f_{t}(\bm{x}_{t}) + c_{t}^{\ast}d(\bm{x}_{t}^{\mathrm{S}},H_{t}^{\ast})\right] 
		\\ \mathrm{s.t.} \quad &\bm{g}_{t}(\bm{x}_{t},H_{t}^{\ast}) \leq \mathbf{0}
		\\ &\bm{x}_{t}^{\mathrm{S}} \geq \mathbf{0}
		\end{split}
		\end{equation}
		\begin{equation}
		\label{3.17}
		c_{t}^{\ast} = \dfrac{-c^{\ast}}{\delta_{t}} - \dfrac{\partial Y(H_{t+1}^{\ast})}{\partial u_{t}}
		\end{equation}
		In (\ref{3.17}), $ c^{\ast} $ is a non-negative constant to be determined by the system conditions in the long run, and $ Y(H_{t+1}) $ is the optimal value of the long-term objective function if the initial SOH is equal to $ H_{t+1} $.
	\end{customthm}
	
	\begin{equation}
	\label{3.20}
	\begin{split}
	J_{t}(c_{t}^{\ast},H_{t}) = &\min\limits_{\bm{x}_{t}} \left[ f_{t}(\bm{x}_{t}) + c_{t}^{\ast}d(\bm{x}_{t}^{\mathrm{S}},H_{t})\right] 
	\\ \mathrm{s.t.} \quad &\bm{g}_{t}(\bm{x}_{t},H_{t}) \leq \mathbf{0}
	\\ &\bm{x}_{t}^{\mathrm{S}} \geq \mathbf{0}
	\end{split}
	\end{equation}
	
	\begin{remark}
		\normalfont
		Theorem \ref{theorem1} says that if the battery operates ($ u_{t}^{\ast} \neq q_{t} $), the optimal dispatch solution is to add a degradation cost term into system objective function with $ c_{t}^{\ast} $ being the optimal time-variant marginal degradation cost. The value of $ c_{t}^{\ast} $ should be set to a discount-factor-adjusted constant plus a term that is dependent on the sensitivity of future profits with respect to the SOH. The value of the constant $ c^{\ast} $, which is essentially an Lagrangian multiplier of Problem B, can be numerically obtained by solving Problem B. $ c_{t} $ is named as the marginal cost of degradation (MCD). 
	\end{remark}
	
	\subsection{Degradation cost evolution and backward algorithm}
	The main purpose of analyzing the optimality conditions of the long-term multi-stage optimization problem in Theorem \ref{theorem1} is to solve the problem and obtain the optimal battery use and system dispatch for each stage. Corollary \ref{corollary1} below presents the backward evolution of MCD, which implies a feasible backward algorithm to solve the problem.
	
	\begin{customcorollary}{1}
		\label{corollary1}
		If Condition \ref{condition1} holds, the optimal solution $ u_{k,t}^{\ast} $ to Problem B (\ref{8})-(\ref{13}) should satisfy either of (\ref{3.21}) and (\ref{3.22}):
		\begin{equation}
		\label{3.21}
		\begin{split}
		u_{k,t}^{\ast} = q_{t} \quad \mathrm{or} \quad u_{k,t+1}^{\ast} = q_{t+1}
		\end{split}
		\end{equation}
		\begin{equation}
		\label{3.22}
		\begin{split}
		c_{t}^{\ast}
		=\dfrac{\delta_{t+1}}{\delta_{t}}
		\left[ 
		c_{t+1}^{\ast} - \dfrac{\partial F_{t+1}(u_{t+1}^{\ast},H_{t+1})}{\partial H_{t+1}}\dfrac{H_{0} - \underline{H} }{ U(H_{0})}
		\right] 
		\end{split}
		\end{equation}
	\end{customcorollary}
	
	\begin{remark}
		\normalfont	
		According to Corollary \ref{corollary1}, the current MCD depends on the future values of MCD and the partial derivative of system cost with respect to battery SOH, and therefore, a backward algorithm is required to find the optimal MCDs, whose details are presented in Algorithm 1. 
	\end{remark}
	
%
 
 

\begin{algorithm}
\DontPrintSemicolon
\SetAlgoLined
\SetKwInOut{Input}{Input}\SetKwInOut{Output}{Output}
\Input{All parameters in $ f_{t}(\boldsymbol{x}_{t}) $, $ \boldsymbol{g}_{t}(\boldsymbol{x}_{t},H_{t}) $, and $ d(\boldsymbol{x}_{t}^{\mathrm{S}},H_{t}) $, which include information of renewable generation and load forecasts, thermal generator capacities and costs, EES capacities, efficiencies, and degradation functions, etc.}
\Output{The optimal degradation costs of EES $ c_{t}^{*} $}
\BlankLine
Choose a reasonable upper bound $ C $ for $ c_{T^{*}}^{*} $ and a required precision $ \Delta c $ \;
\ForEach{$ T' \in [1,N] $}{
Initialize $ c_{T'} = 0 $ \; 
\While{$ c_{T'}<C $}{
    Initialize $ t = T' + 1 $ and $ H_{T'+1} = \underline{H} $\;
    \While{$ H_{t} < H_{0}$  and  $t \geq 1 $}{
    	$ t = t - 1 $ \;
	    Solve Problem C to get the optimal solution $ x'_{t} $ \;
	    $ u'_{t} \leftarrow d(x'_{t}, H_{t+1}) $ \;
	    $ H_{t} \leftarrow H_{t+1} +\dfrac{H_{0}-\underline{H}}{U(H_{0})} u'_{t} $ \;
	    Calculate $ F_{t}(u'_{t},H_{t})$ using Problem A \;	     
	    Calculate $ \dfrac{\partial{F_{t}(u'_{t},H_{t})}}{\partial{H_{t}}}$ \; 
	    $  c_{t-1}(T',c_{T'}) \leftarrow \dfrac{\delta_{t}}{\delta_{t-1}}\left[c_{t}(T',c_{T'}) - \dfrac{\partial{F_{t}(u'_{t},H_{t})}}{\partial{H_{t}}} \dfrac{H_{0}-\underline{H}}{U(H_{0})} \right]   $ \; 
	}
	$ t_{0} \leftarrow t $ \;   
	\If{$ t_{0} = 1 $}{
		$ y(T',c_{T'}) \leftarrow \sum\limits_{ t \in [t_{0},T'] } \delta_{t}F_{t}(u'_{t},H_{t})$
		}

	$ c_{T'} \leftarrow c_{T'} + \Delta c $ 
} 
$ T^{*},c_{T^{*}}^{*} \leftarrow \argmin\limits_{ T',c_{T'} }  y(T',c_{T'}) $ \; 
$ c_{t}^{*} \leftarrow c_{t}(T^{*},c_{T^{*}}^{*}) $

}
\caption{Backward Algorithm}
\end{algorithm}


	Algorithm 1 starts from the last year of the battery life (which is also a variable to be optimized), sequentially calculates the MCD, battery usage and SOH, system cost and its partial derivative with respect to battery SOH, and ends up with the first year of battery life. The MCD series with the minimum system cost are the optimal ones. 
	
	Algorithm 1 may be time-consuming when the system is large. However, the long-term problem only needs to be solved when there are significant unforeseen changes in the system resources or topology, which might affect the optimal values of MCD. In the day-ahead or real-time dispatch, only Problem C needs to be solved with implementing the optimal MCD. Besides, algorithms that are more efficient can be designed if more properties of $ F_{t}(u_{t},H_{t}) $ are known or previous optimal solutions are available for warm starts.	
	
	\section{Case Study on Economic Dispatch}
	\label{sec:studytwo:case_study}
	In this section, an optimization model and simulation results for the economic dispatch of a sample power system are presented. The sample power system consists of one wind farm, one thermal power plant, loads, and one battery system. In the simulation results, the evolution of the optimal MCDs is demonstrated, and the life-cycle cost saving using the optimal MCD series is compared with that using a degradation cost derived from battery capital cost. The effect of incorporating the SOH term (the second term in (\ref{3.16}) and (\ref{3.22})) in the optimal MCD is also examined through life-cycle cost saving comparison. 
	
	\subsection{Model formulations}
	Here the formulations of Problem A and C for the case study are presented, respectively. Problem A includes (\ref{3.23})-(\ref{3.37}), while Problem C includes  (\ref{3.24})-(\ref{3.36}) and (\ref{3.38})-(\ref{3.39}).
	
	\subsubsection{ Degradation-constrained model (Problem A)}
	The economic dispatch model is a cost-minimizing problem, with the power output schedules of thermal power plants, controllable loads, and battery units as its decision variables, as equation (\ref{3.23}). The time horizon of each dispatch decision is usually a day, represented by $ \Delta t $, and $ t $ is the day index.
	\begin{equation}
	\label{3.23}
	\begin{split}
	F_{t}(u_{k,t},H_{k,t}) &= \min\limits_{x^{\mathrm{G}}_{m,h},x^{\mathrm{L}}_{n,h},x^{\mathrm{S+}}_{k,h},x^{\mathrm{S-}}_{k,h}} f_{t}(\bm{x}_{t})
	\\ &= \min\limits_{x^{\mathrm{G}}_{m,h},x^{\mathrm{L}}_{n,h},x^{\mathrm{S+}}_{k,h},x^{\mathrm{S-}}_{k,h}}
	\left(\sum\limits_{m \in \Omega_{m}} C^{\mathrm{G}}_{m}
	+ \sum\limits_{n \in \Omega_{n}} C^{\mathrm{L}}_{n}
	\right) 
	\end{split}
	\end{equation}
	
	The total system cost during day $ t $ is the sum of the costs of all thermal power plants, controllable loads. The cost of the $ m $ th thermal power plant, $ C^{\mathrm{G}}_{m} $, is assumed to be a quadratic function of its generation output at hour $ h $ denoted by $ x^{\mathrm{G}}_{m,h} $, as in equation (\ref{3.24}). $ a^{\mathrm{G}}_{m} $ and $ b^{\mathrm{G}}_{m} $ are the coefficients of the cost function of the $ m $ th thermal power plant, and $ \Omega_{m} $ is the set of thermal power plants in the system, whose cardinality is 1 in this case study.
	\begin{equation}
	\label{3.24}
	\begin{split}
	C^{\mathrm{G}}_{m} &= \sum\limits_{h \in [t,t+\Delta t]}
	\left[ a^{\mathrm{G}}_{m}\left(x^{\mathrm{G}}_{m,h} \right)^{2} +
	b^{\mathrm{G}}_{m} x^{\mathrm{G}}_{m,h}  \right] 
	\qquad \forall m \in \Omega_{m}
	\end{split}
	\end{equation}
	
	The cost of the $ n $ th controllable load, $ C^{\mathrm{L}}_{n} $, is also assumed to be a quadratic function of its load reduction at hour $ h $ denoted by $ x^{\mathrm{L}}_{n,h} $, as in equation (\ref{3.25}). $ a^{\mathrm{L}}_{n} $ and $ b^{\mathrm{L}}_{n} $ are the coefficients of the cost function of the $ n $ th controllable load, and $ \Omega_{n} $ is the set of controllable loads in the system, whose cardinality is 1 in this case study.
	\begin{equation}
	\label{3.25}
	\begin{split}
	C^{\mathrm{L}}_{n} &= \sum\limits_{h \in [t,t+\Delta t]}
	\left[ a^{\mathrm{L}}_{n}\left(x^{\mathrm{L}}_{n,h} \right)^{2} +
	b^{\mathrm{L}}_{n} x^{\mathrm{L}}_{n,h}  \right] 
	\qquad \forall n \in \Omega_{n}
	\end{split}
	\end{equation}
	
	The power balance constraint of the system is as equation (\ref{3.26}). In equation (\ref{3.26}), $ x^{\mathrm{S-}}_{k,h} $ and $ x^{\mathrm{S+}}_{k,h} $ are charging and discharging power of the $ k $ th battery ; $ W_{\omega,h} $ is the forecasted available wind energy of the $ \omega $ th wind farm during time $ h $; $ L_{n,h} $ is the forecasted load level at node $ n $; $ \Omega_{\rho} $ is the set of wind farms in the system with a cardinality of 1 in this case.
	\begin{equation}
	\label{3.26}
	\begin{split}
	&\left[ \sum\limits_{m \in \Omega_{m}} x^{\mathrm{G}}_{m,h} + \sum\limits_{k \in \Omega_{k}} \left( x^{\mathrm{S+}}_{k,h} - x^{\mathrm{S-}}_{k,h}\right)  \right] \Delta h + \sum\limits_{\omega \in \Omega_{\omega}} W_{\omega,h}
	\\\geq &
	\sum\limits_{n \in \Omega_{n}} \left(L_{n,h} - x^{\mathrm{L}}_{n,h} \right) \Delta h 
	\qquad \forall h \in [t,t+\Delta t]
	\end{split}
	\end{equation}
	
	The control variables have physical upper limits as equations (\ref{3.27}) to (\ref{3.29}), where $ \overline{x}^{\mathrm{G}}_{m} $, $ \overline{x}^{\mathrm{L}}_{n} $, and $ \overline{x}^{\mathrm{S}}_{k,t} $ represent the limits of thermal generation plants, controllable loads, and battery units, respectively.
	\begin{equation}
	\label{3.27}
	\begin{split}
	0 \leq x^{\mathrm{G}}_{m,h}\leq \overline{x}^{\mathrm{G}}_{m}
	\qquad \forall m \in \Omega_{m}, h \in [t,t+\Delta t]
	\end{split}
	\end{equation}
	\begin{equation}
	\label{3.28}
	\begin{split}
	0 \leq x^{\mathrm{L}}_{n,h}\leq \overline{x}^{\mathrm{L}}_{n}
	\qquad \forall n \in \Omega_{n}, h \in [t,t+\Delta t]
	\end{split}
	\end{equation}
	\begin{equation}
	\label{3.29}
	\begin{split}
	0 \leq x^{\mathrm{S+}}_{k,h},x^{\mathrm{S-}}_{k,h}\leq \overline{x}^{\mathrm{S}}_{k,t}
	\qquad \forall k \in \Omega_{k}, h \in [t,t+\Delta t]
	\end{split}
	\end{equation}
	
	The energy constraints of battery are modelled in equations (\ref{3.30}) to (\ref{3.32}). The energy level of battery at time $ h+1 $, $ e^{\mathrm{S}}_{k,h+1} $, is expressed as a function of the energy level at time $ h $, $ e^{\mathrm{S}}_{k,h} $, and the charging/discharging output at time $ h $, as equation (\ref{3.30}), where $ \rho $ is the self-discharge rate, and $ \eta_{k,t} $ is the charge/discharge efficiency.
	\begin{equation}
	\label{3.30}
	\begin{split}
	e^{\mathrm{S}}_{k,h+1} = (1-\rho)e^{\mathrm{S}}_{k,h} - \dfrac{x^{\mathrm{S+}}_{k,h}}{\eta_{k,t}}\Delta h + x^{\mathrm{S-}}_{k,h}\eta_{k,t}\Delta h
	\\ \forall k \in \Omega_{k}, h \in [t,t+\Delta t]
	\end{split}
	\end{equation}
	
	The energy level of the battery has to be kept within its capacity, as equation (\ref{3.31}), where $ \overline{e}^{\mathrm{S}}_{k,t} $ is the energy capacity of the $ k $ th battery unit.
	\begin{equation}
	\label{3.31}
	\begin{split}
	0 \leq e^{\mathrm{S}}_{k,h}\leq \overline{e}^{\mathrm{S}}_{k,t}
	\qquad \forall k \in \Omega_{k}, h \in [t,t+\Delta t]
	\end{split}
	\end{equation}
	
	To reserve flexibility, the initial and the final energy levels are set to be equal for the dispatch horizon $ \left[ t,t+\Delta t \right]  $, as equation (\ref{3.32}).
	\begin{equation}
	\label{3.32}
	\begin{split}
	e^{\mathrm{S}}_{k,t} = e^{\mathrm{S}}_{k,t+\Delta t}
	\end{split}
	\end{equation}
	
	The energy capacity $ \overline{e}^{\mathrm{S}}_{k,t} $ and the internal impedance $ Z_{k,t} $ of battery are functions of the SOH $ H_{k,t} $, while the power capacity $ \overline{x}^{\mathrm{S}}_{k,t} $, and the efficiency $ \eta_{k,t} $ of battery are functions of the internal impedance, as equations (\ref{3.33})-(\ref{3.36}) \cite{shi18-2nQfj}. $ \overline{e}^{\mathrm{S}}_{k,0} $, $ Z_{k,0} $, $ H_{k,0} $, and $ \overline{x}^{\mathrm{S}}_{k,t} $ are the initial energy capacity, internal impedance, SOH, and power capacity, respectively. $ \overline{Z}_{k} $ is the internal impedance at the end of battery life. 
	\begin{equation}
	\label{3.33}
	\begin{split}
	\overline{e}^{\mathrm{S}}_{k,t} = H_{k,t}\overline{e}^{\mathrm{S}}_{k,0}
	\end{split}
	\end{equation}
	\begin{equation}
	\label{3.34}
	\begin{split}
	Z_{k,t} = Z_{k,0} + (\overline{Z}_{k} - Z_{k,0})\dfrac{H_{k,0}-H_{k,t}}{H_{k,0}-\underline{H}_{k}}
	\end{split}
	\end{equation}
	\begin{equation}
	\label{3.35}
	\begin{split}
	\overline{x}^{\mathrm{S}}_{k,t} = \dfrac{Z_{k,0}}{Z_{k,t}}\overline{x}^{\mathrm{S}}_{k,0}
	\end{split}
	\end{equation}
	\begin{equation}
	\begin{split}
	\label{3.36}
	\eta_{k,t} &= \dfrac{1}{1+\dfrac{Z_{k,t}}{Z_{k,0}}\dfrac{1-\eta_{k,0}}{\eta_{k,0}}}
	\end{split}
	\end{equation}
	
	For degradation-constrained model, the battery usage/degradation is capped, as equation (\ref{3.37}). The battery degradation/usage is assumed to be linearly dependent on its charging and discharging power output, $ x^{\mathrm{S-}}_{k,h} $ and $ x^{\mathrm{S+}}_{k,h} $. The degradation can be a more complicated convex non-linear function or piece-wise linear function of the charging and discharging schedules to fit the real degradation characteristics of a specific battery system.
	\begin{equation}
	\begin{split}
	\label{3.37}
	\sum\limits_{h \in [t,t + \Delta t]}\left( x^{\mathrm{S+}}_{k,h} + x^{\mathrm{S-}}_{k,h} \right) \Delta h + q_{k,t} \leq u_{k,t}
	\end{split}
	\end{equation}
	
	\subsubsection{ Degradation-cost model (Problem C)}
	For degradation-cost model, a battery cost term is added to the objective function, as equation (\ref{3.38}).
	\begin{equation}
	\label{3.38}
	\begin{split}
	&J_{t}(c_{k,t},H_{k,t}) 
	\\=&  \min\limits_{x^{\mathrm{G}}_{m,h},x^{\mathrm{L}}_{n,h},x^{\mathrm{S+}}_{k,h},x^{\mathrm{S-}}_{k,h}} \left[ f_{t}(\bm{x}_{t} ) + \sum\limits_{k \in \Omega_{k}} C^{\mathrm{S}}_{k}\right] 
	\\ =& \min\limits_{x^{\mathrm{G}}_{m,h},x^{\mathrm{L}}_{n,h},x^{\mathrm{S+}}_{k,h},x^{\mathrm{S-}}_{k,h}}
	\left(\sum\limits_{m \in \Omega_{m}} C^{\mathrm{G}}_{m}
	+ \sum\limits_{n \in \Omega_{n}} C^{\mathrm{L}}_{n}
	+ \sum\limits_{k \in \Omega_{k}} C^{\mathrm{S}}_{k}
	\right) 
	\end{split}
	\end{equation}
	
	The degradation cost of the $ k $ th battery unit $ C_{k}^{\mathrm{S}} $ is proportional to its degradation/usage, as equation (\ref{3.39}).
	\begin{equation}
	\begin{split}
	\label{3.39}
	C_{k}^{\mathrm{S}} = c_{k,t}\left[ \sum\limits_{h \in [t,t + \Delta t]}\left( x^{\mathrm{S+}}_{k,h} + x^{\mathrm{S-}}_{k,h} \right) \Delta h + q_{k,t} \right] 
	\\ \forall k \in \Omega_{k}
	\end{split}
	\end{equation}
	
	Equations (\ref{3.24})-(\ref{3.36}) also belong to this model.
	
	\subsection{Data}
\begin{figure}[ht]
		\centering
		\includegraphics[width=\textwidth/2]{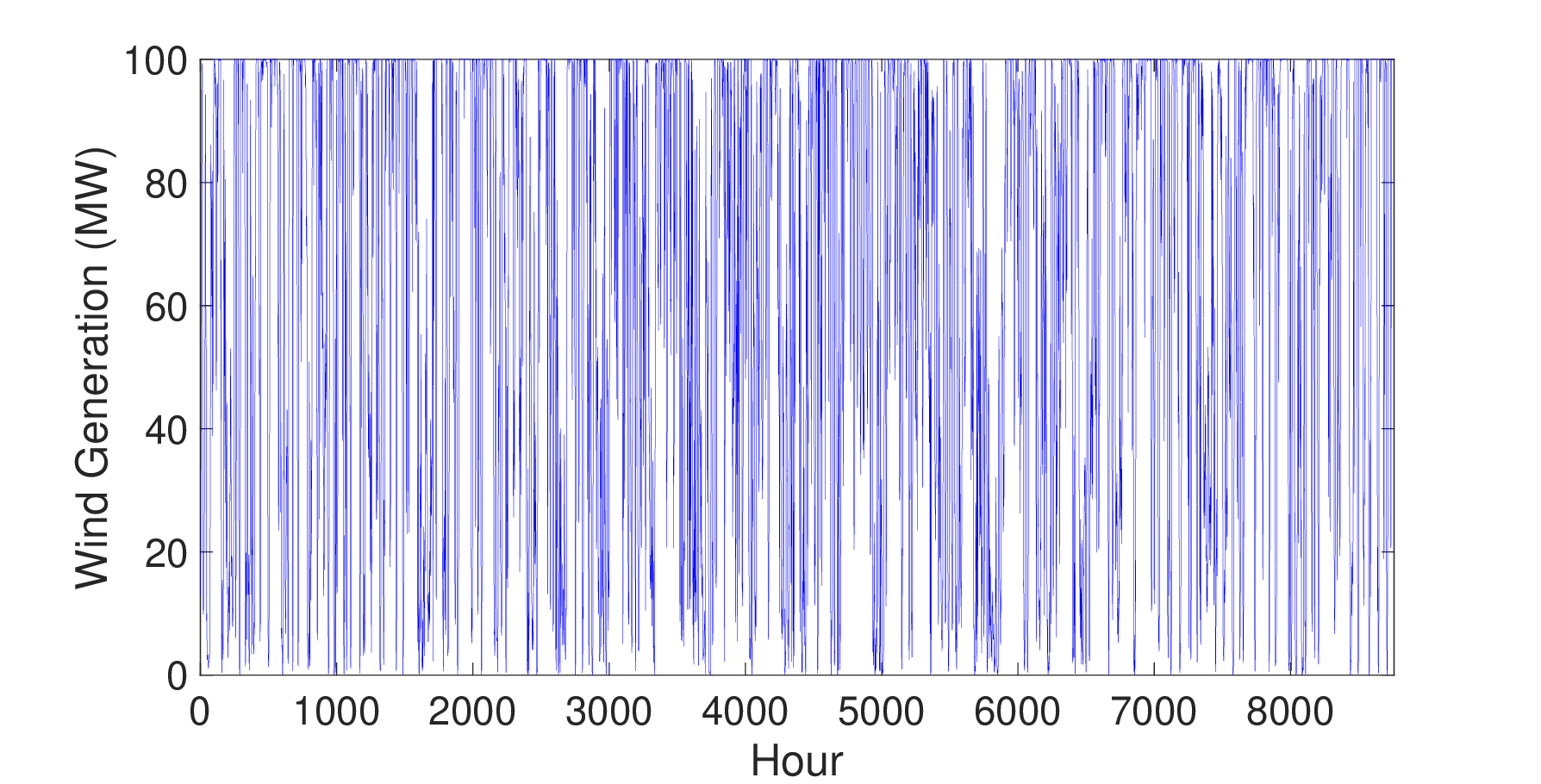}
		\caption{Annual wind generation profile based on a selected location in Texas in 2015.}
		\label{fig-ch3-1}
	\end{figure}

\begin{figure}[ht]
		\centering
		\includegraphics[width=\textwidth/2]{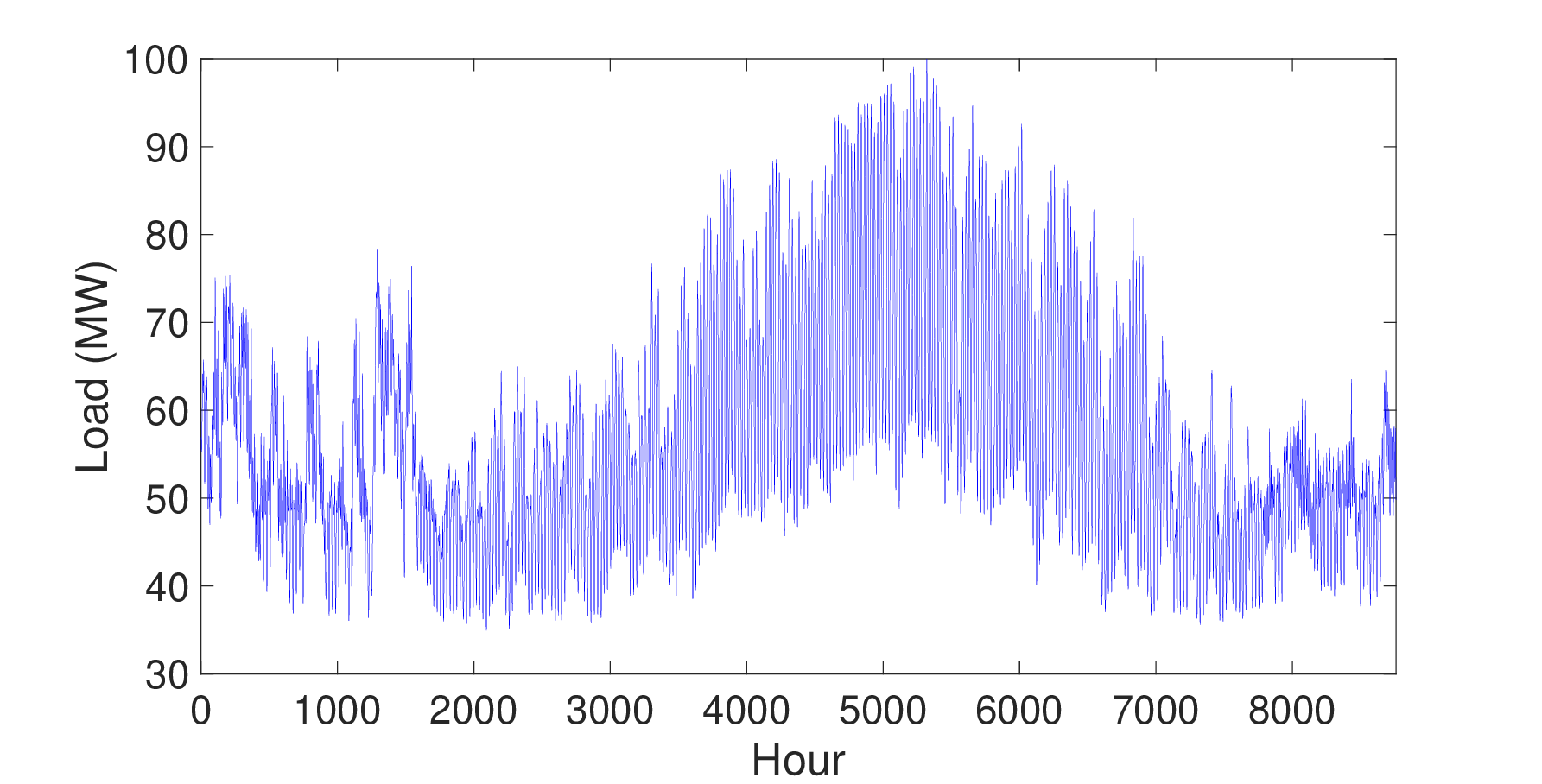}
		\caption{Annual load profile rescaled from the annual total load profile in ERCOT in 2015.}
		\label{fig-ch3-2}
	\end{figure}

	The wind generation profile is presented in Fig. \ref{fig-ch3-1}, which is produced based on the historical 100-m wind speed data in 2015 in a selected location in Texas. The wind capacity factor is approximately 63\%. The load profile of the system is presented in Fig. \ref{fig-ch3-2}, which is scaled from the annual total load profile in ERCOT in 2015. The average load is approximately 57 MW. Both the wind and load profiles are assumed to be the same for each year across the lifetime of the battery. The parameters of the thermal power plant, controllable load, and battery are summarized in Table \ref{tab-ch3-1} and \ref{tab-ch3-2}. The life-cycle usage limit of the battery, measured by energy throughput, is 1.2 TWh, which is 3000 full cycles for a 200-MWh lithium-ion battery \cite{sarasketa-zabala15,wu18}. After 3000 full cycles, the remaining energy capacity of the battery decreases to 70\% of the original capacity and the impedance increases to 200\% of the original [38]. The discount rate is set to 7\%.
	\begin{table}[ht]
		\renewcommand{\arraystretch}{1.3}
		\caption{Parameters of the thermal power plant and controllable load in the system.}
		\label{tab-ch3-1}
		\centering
		\begin{tabular}{p{0.4in}<{\centering}p{0.4in}<{\centering}p{0.4in}<{\centering}p{0.4in}<{\centering}p{0.4in}<{\centering}p{0.4in}<{\centering}}
			\hline\hline
			$ \overline{x}^{\mathrm{G}} $ & 
			$ a^{\mathrm{G}} $&
			$ b^{\mathrm{G}} $&
			$ \overline{x}^{\mathrm{L}} $&
			$ a^{\mathrm{L}} $&
			$ b^{\mathrm{L}} $\\
			\midrule
			100 MW &
			0.1 \$/MW$ ^{2} $ &
			30 \$/MW &
			10 MW &
			0.1 \$/MW$ ^{2} $ &
			70 \$/MW\\		
			\hline\hline
		\end{tabular}
	\end{table}
	\begin{table}[ht]
		\renewcommand{\arraystretch}{1.3}
		\caption{Parameters of the battery.}
		\label{tab-ch3-2}
		\centering
		\begin{tabular}{p{0.4in}<{\centering}p{0.4in}<{\centering}p{0.4in}<{\centering}p{0.4in}<{\centering}p{0.4in}<{\centering}p{0.4in}<{\centering}}
			\hline\hline
			$ \overline{x}^{\mathrm{S}}_{0} $ & 
			$ a^{\mathrm{S}}_{0} $&
			$ \eta_{0} $&
			$ \sigma $&
			$ U $&
			$ q $\\
			\midrule
			50 MW &
			200 MWh &
			95\% &
			0 &
			1.2 TWh &
			50 MWh\\		
			\hline\hline
		\end{tabular}
	\end{table}
	
	\subsection{Results}
\begin{figure}
		\centering
		\includegraphics[width=\textwidth/2]{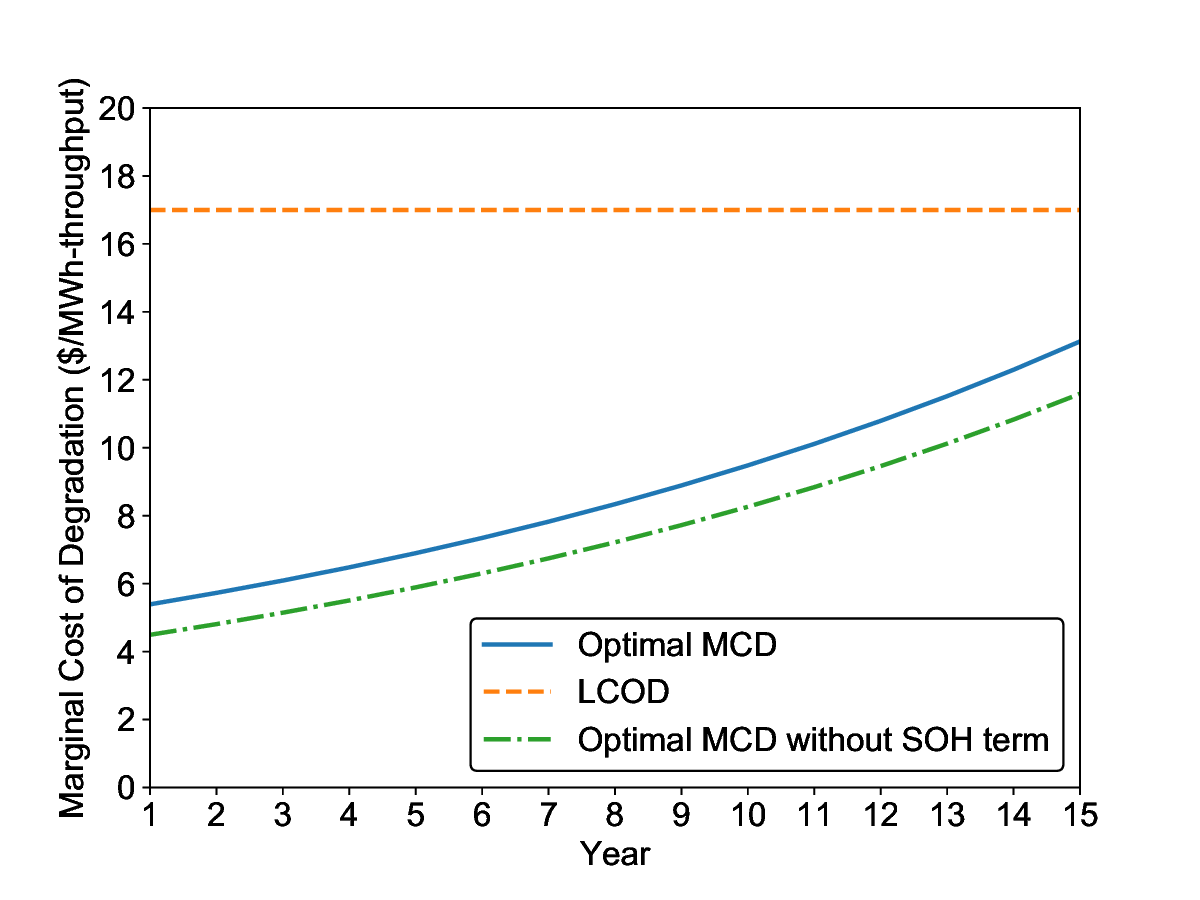}
		\caption{The optimal MCDs with and without the term concerning SOH and the LCOD.}
		\label{fig-ch3-3}
	\end{figure}

	By implementing Algorithm 1, the operation of the sample system is simulated, and the optimal MCDs as well as the total contribution of the battery in terms of cost saving are calculated. The optimal MCDs, or say, the optimal unit degradation costs, from Year 1 to Year 15 are presented in Fig. \ref{fig-ch3-3} (blue solid line). It can be observed that the optimal MCD increases from approximately \$7/MWh-throughput to \$13/MWh-throughput. This evolving trend implies that the battery should be used more frequently and intensively in the early years of its lifetime, if the microgrid configuration and generating costs do not change. The total system cost saving from battery operation is approximately \$9.4 million.
	
	The levelized cost of degradation (LCOD) method im-plemented in many existing studies, uses an amortized proportion of initial capital cost or future replacement cost as the unit degradation cost. Assuming that the unit-capacity capital cost of the lithium-ion battery system is \$200/kWh, that the degradation is uniformly allocated over 15 years, and that the ratio of total depreciation to capital cost is 30\%, the LCOD is approximately \$17/MWh, as plotted in Fig. \ref{fig-ch3-3} (orange dash line). The total system cost saving from storage operation using this time-invariant unit degradation cost is \$4.1 million, which is only 44\% of the proposed method. The significant loss in cost saving is because the LCOD method does not maximize the life-cycle revenue but requires the battery to operate only when the potential marginal benefit is higher than the average unit degradation cost, the LCOD, and to wait for a better benefit opportunity otherwise. The LCOD method does not consider the fact that the battery has a calendar life\textemdash some less profitable opportunities should also be captured before the life ends.
	
	The optimal MCD series without incorporating the second SOH term $ - \dfrac{\partial Y(H_{k,t+1})}{\partial u_{k,t}} $ in (\ref{3.17}) is also plotted in Fig. \ref{fig-ch3-3}, which does not deviate too much from the optimal MCD series with the SOH term. The MCD increases a bit after considering the effect of current usage on future SOH, because the battery tends to incur less degradation at current to increase the future SOH. The total system cost saving without the SOH term is \$9.1 million, which is 3\% less than that with the SOH term.
	
\begin{figure}
		\centering
		\includegraphics[width=\textwidth/2]{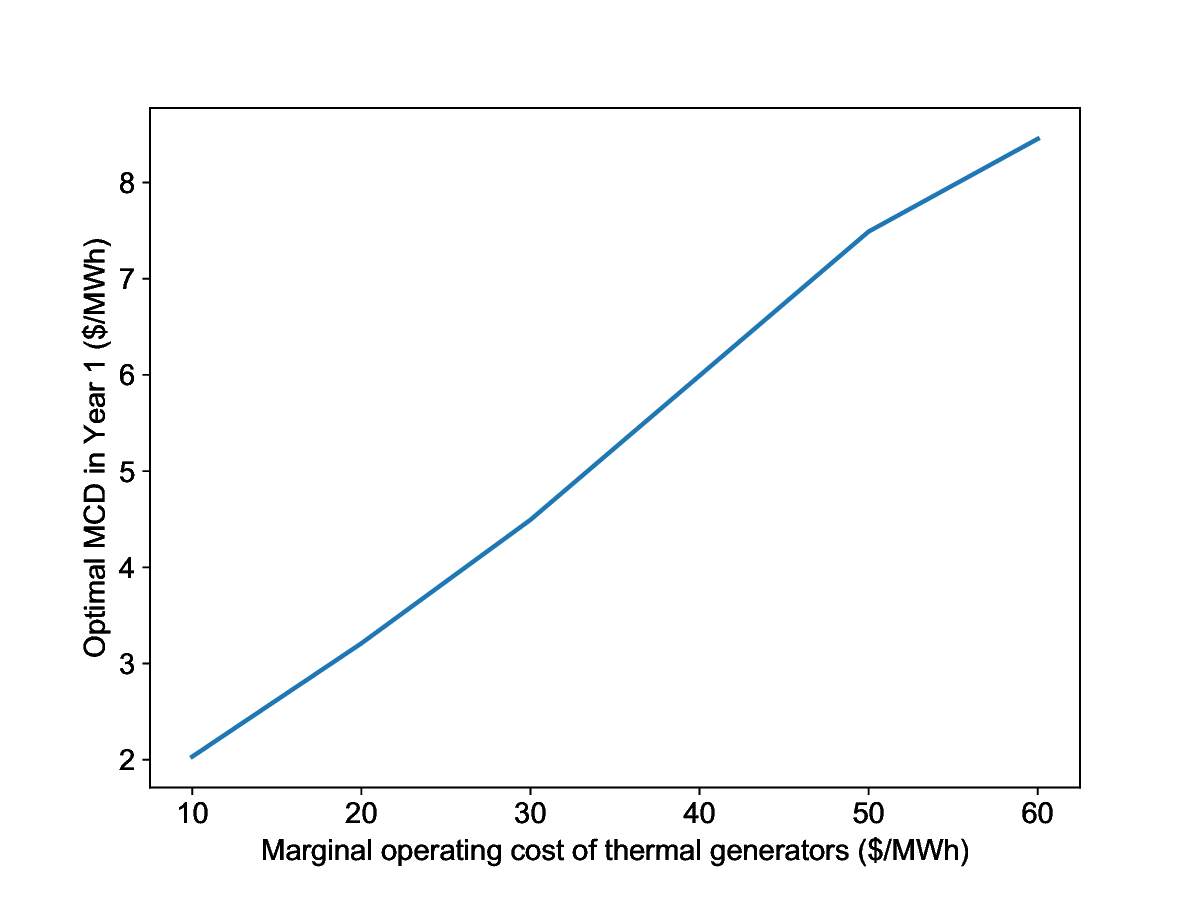}
		\caption{Different optimal marginal costs of degradation in the first year as the marginal cost of the thermal power plant varies.}
		\label{fig-ch3-4}
	\end{figure}

	To find out how the degradation cost of the battery should vary as the system parameters change, the system operation is simulated with different values of the marginal cost of the thermal power plant, $ b^{\mathrm{G}} $. As observed in Fig. \ref{fig-ch3-4}, the optimal MCD is linearly increasing as the marginal cost of the thermal power plant increases. This indicates that the degradation cost of battery should be positively dependent on the marginal cost of the system\textemdash if the marginal cost of the system changes, the MCD should also change. 
	
	\section{Conclusions}
	\label{sec:studytwo:conclusions}
	Dispatch models are the center of smart grid operation, and better dispatch models could imply significant cost savings. This study derives the intertemporal condition and the degradation cost for a power system with battery, with mathematical proofs of their optimality. It is proved that there exists a value, named MCD, that should be adjusted over time and implemented into the cost function of battery, as if it is the variable unit operational cost incurred by degradation. The MCD is independent of the capital cost of the battery at the operational stage. In case studies, the optimal outcome of the proposed method is compared with the LCOD method used in existing studies as well as with the approach without considering the SOH-related term. It turns out that the proposed method can generate much more benefits than the LCOD method. The benefit improvement of incorporating the SOH term is only marginal in this case, which implies that the SOH term may be ignored in some cases to reduce the computational complexity brought by it.
	
	The cost function of battery derived in this paper could serve as a foundation for dispatch and planning problems concerning battery.

	\section*{Appendix}
	\subsection{Proof of Theorem \ref{theorem1}}
	Because our main results are only concerned with the optimality conditions for battery storage within its lifetime, when $ \beta_{t}=1 $, Problem B can be transformed to the following, with the SOH evolution embedded in the objective function:
	\begin{equation}
	\label{e41}
	\begin{split}
		Y(H_{0})  = &\min\limits_{u_{t}} \sum\limits_{t \leq N} \delta_{t} F_{t}\left(u_{t},H_{t}\left(H_{0},\sum\limits_{\tau < t}u_{\tau} \right)  \right)  	\\ 
		\mathrm{s.t.} \qquad &\sum\limits_{t \leq N} u_{t} \leq U(H_{0})\\ 
		&u_{t} \geq q_{t} \qquad\qquad\qquad\qquad \forall t \leq N
	\end{split}
	\end{equation}
	Let's first introduce two Lemmas:
	\begin{customlemma}{1}
		\label{lemma1}
		If $F_{t}\left(u_{t},H_{t} \right)$ is differentiable with respect to $ u_{t} $, the optimal solution $u_{t}^{\ast}$ to (\ref{e41}) should satisfy either of (\ref{A.1}) and (\ref{A.2}):
		\begin{align}
		u_{t}^{\ast} = q_{t}
		\label{A.1}
		\end{align}
		\begin{align}
		\dfrac{\partial F_{t}\left(u_{t}^{\ast},H_{t} \right) }{\partial u_{t}} = -c_{t}^{\ast}
		=\dfrac{-c^{\ast}}{\delta_{t}} - \dfrac{\partial Y(H_{t+1})}{\partial u_{t}}
		\label{A.2}
		\end{align}
	\end{customlemma}
	
	\begin{proof}
		$ \forall t \leq N $, (\ref{e41}) can be reformed as:
		\begin{equation}
		\begin{split}
		\label{3.40}
		Y(H_{0})  = &\min\limits_{u_{t}} \sum\limits_{\tau \leq t} \delta_{\tau} F_{\tau}\left(u_{\tau},H_{\tau}\right) + \delta_{t}Y(H_{t+1})
		\\ \mathrm{s.t.} \quad &\sum\limits_{\tau \leq t}  u_{\tau}
		\leq U(H_{0}) - U(H_{t+1}) : \qquad c \geq 0
		\\ &u_{\tau} \geq q_{\tau}: \qquad \alpha_{\tau} \geq 0 \qquad\qquad\qquad \forall \tau \leq t
		\end{split}
		\end{equation}
		
		The Lagrangian function of problem (\ref{3.40}) is:
		\begin{equation}
		\begin{split}
		\label{3.41}
		L = &\sum\limits_{\tau \leq t} \delta_{\tau} F_{\tau}\left(u_{\tau},H_{\tau}\right) + \delta_{t}Y(H_{t+1})
		\\+& c\left[ \sum\limits_{\tau \leq t}  u_{\tau} - U(H_{0}) + U(H_{t+1}) \right] + \sum\limits_{\tau \leq t} \alpha_{\tau}(q_{\tau} - u_{\tau})
		\end{split}
		\end{equation}
		where $ c $ and $ \alpha_{t} $ are Lagrangian multipliers. Then, three of the first-order Karush–Kuhn–Tucker (KKT) conditions of problem (\ref{3.40}) are equations (\ref{3.42})-(\ref{3.44}):
		\begin{equation}
		\begin{split}
		\label{3.42}
		\dfrac{\partial L}{\partial u_{t}} &= 
		\delta_{t} \dfrac{\partial F_{t}\left(u_{t}^{\ast},H_{t}^{\ast} \right)}{\partial u_{t}}
		+ \delta_{t} \dfrac{\partial Y(H_{t+1}^{\ast})}{\partial u_{t}}
		+ c^{\ast} - \alpha_{t}^{\ast}
		\\&=0
		\end{split}
		\end{equation}
		\begin{equation}
		\begin{split}
		\label{3.43}
		\alpha_{t}^{\ast} (q_{t} - u_{t}^{\ast}) = 0
		\end{split}
		\end{equation}
		\begin{equation}
		\begin{split}
		\label{3.44}
		u_{t}^{\ast} \geq q_{t}
		\end{split}
		\end{equation}
		
		If $ u_{t}^{\ast} > q_{t} $, according to equation (\ref{3.43}), the following holds:
		\begin{equation}
		\begin{split}
		\label{3.45}
		\alpha_{t} = 0
		\end{split}
		\end{equation}
		
		Then substituting equation (\ref{3.45}) into equation (\ref{3.42}), (\ref{A.2}) is obtained.
	\end{proof}
	
	\begin{customlemma}{2}
		\label{lemma2}
		If $F_{t}\left(u_{t},H_{t} \right)$ is differentiable and convex with respect to $ u_{t} $, and $ u_{t}^{\ast} $ and $ c_{t}^{\ast} $ are the optimal solutions to Problem B, then Problem A' and Problem C are equivalent (having identical optimal solution sets).
		
		Problem A':
		\begin{equation}
		\label{3.19}
		\begin{split}
		F_{t}(u_{t}^{\ast},H_{t}^{\ast}) = &\min\limits_{\bm{x}_{t}} f_{t}(\bm{x}_{t})
		\\ \mathrm{s.t.} \quad &\bm{g}_{t}(\bm{x}_{t},H_{t}^{\ast}) \leq \mathbf{0}
		\\ &d(\bm{x}_{t}^{\mathrm{S}},H_{t}^{\ast}) \leq u_{t}^{\ast}
		\\ &\bm{x}_{t}^{\mathrm{S}} \geq \mathbf{0}
		\end{split}
		\end{equation}
	\end{customlemma}
	
	\begin{proof}
		First, Problem D is introduced as:
		\begin{equation}
		\label{3.46}
		\begin{split}
		&\min\limits_{u_{t}}\min\limits_{\bm{x}_{t}} \left[ f_{t}(\bm{x}_{t}) + c_{t}^{\ast}u_{t}\right] 
		\\ \mathrm{s.t.} \quad &\bm{g}_{t}(\bm{x}_{t},H_{t}^{\ast}) \leq \mathbf{0}
		\\ &d_{t}(\bm{x}_{t}^{\mathrm{S}},H_{t}^{\ast}) \leq u_{t}
		\\ &\bm{x}_{t}^{\mathrm{S}} \geq \mathbf{0}
		\end{split}
		\end{equation}
		
		Problem D can be reformulated as below:
		
		Problem D:
		\begin{equation}
		\label{3.47}
		\begin{split}
		&\min\limits_{u_{t}} D_{t}(u_{t},H_{t}^{\ast})
		\end{split}
		\end{equation}
		
		Sub-Problem D:
		\begin{equation}
		\label{3.48}
		\begin{split}
		&D_{t}(u_{t},H_{t}^{\ast}) = \min\limits_{\bm{x}_{t} \in \mathbf{\Phi}} \left[ f_{t}(\bm{x}_{t}) + c_{t}^{\ast}u_{t}\right] 
		\\ &\mathbf{\Phi} =  \left\lbrace \bm{x}_{t}|\bm{g}_{t}(\bm{x}_{t},H_{t}^{\ast}) \leq \mathbf{0},d_{t}(\bm{x}_{t}^{\mathrm{S}},H_{t}^{\ast}) \leq u_{t},\bm{x}_{t}^{\mathrm{S}} \geq \mathbf{0}\right\rbrace 
		\end{split}
		\end{equation}
		where $ D_{t}(u_{t},H_{t}^{\ast}) $ is the optimal value of sub-problem D. According to (\ref{12}) and (\ref{3.48}), $ u_{t}^{\ast } \geq q_{t} $ is also an implied constraint.
		
		It can be observed that $ d(\bm{x}_{t}^{\mathrm{S}},H_{t}^{\ast}) = u_{t} $ is always a necessary optimal condition to Problem D, as $ D_{t}(d(\bm{x}_{t}^{\mathrm{S}},H_{t}^{\ast}),H_{t}^{\ast}) \leq D_{t}(u_{t},H_{t}^{\ast}) $. Therefore, the optimal solution set to Problem C is the same as the optimal solution set to Problem D, considering that Problem C and D have the same feasible region and objective function. Each optimal solution to problem C is a feasible solution to Problem D, and vice versa. 
		
		$ D_{t}(u_{t},H_{t}^{\ast}) $ has the following relation with $ F_{t}(u_{t},H_{t}^{\ast}) $:
		\begin{equation}
		\label{3.49}
		\begin{split}
		&D_{t}(u_{t},H_{t}^{\ast}) = F_{t}(u_{t},H_{t}^{\ast}) + c_{t}^{\ast}u_{t}
		\end{split}
		\end{equation}
		
		The property expressed in equation (\ref{3.49}) is obvious because when $ u_{t} $ is fixed, $ c_{t}^{\ast}u_{t} $ is a constant, and $ \min\limits_{\bm{x}_{t} \in \mathbf{\Phi}} \left[ f_{t}(\bm{x}_{t}) + c_{t}^{\ast}u_{t}\right]  $ has the same optimal solution set with $ \min\limits_{\bm{x}_{t} \in \mathbf{\Phi}} \left[ f_{t}(\bm{x}_{t}) \right]  $. Implementing equation (\ref{3.49}) into Problem D:
		\begin{equation}
		\label{3.50}
		\begin{split}
		&\min\limits_{u_{t}} F_{t}(u_{t},H_{t}^{\ast}) + c_{t}^{\ast}u_{t}
		\\&\mathrm{s.t.} \quad u_{t} \geq q_{t}
		\end{split}
		\end{equation}
		
		The KKT conditions of Problem D are equations (\ref{3.51})-(\ref{3.54}):
		\begin{equation}
		\begin{split}
		\label{3.51}
		\dfrac{\partial L_{\mathrm{D}}(u_{t}^{\ast D})}{\partial u_{t}} =&
		\dfrac{\partial \left[ F_{t}\left(u_{t}^{\ast D},H_{t}^{\ast}\right) + c_{t}^{\ast}u_{t} + \nu_{t}^{\ast}\left( q_{t} - u_{t}^{\ast D} \right)   \right] }{\partial u_{t}} 
		\\=& 0
		\\ \Leftrightarrow \quad \quad &\dfrac{\partial  F_{t}\left(u_{t}^{\ast D},H_{t}^{\ast}\right) }{\partial u_{t}} = -c_{t}^{\ast} + \nu_{t}^{\ast}
		\end{split}
		\end{equation}
		\begin{equation}
		\begin{split}
		\label{3.52}
		\nu_{t}^{\ast} (q_{t} - u_{t}^{\ast D}) = 0
		\end{split}
		\end{equation}
		\begin{equation}
		\begin{split}
		\label{3.53}
		u_{t}^{\ast D} \geq q_{t}
		\end{split}
		\end{equation}
		\begin{equation}
		\begin{split}
		\label{3.54}
		\nu_{t}^{\ast} \geq 0
		\end{split}
		\end{equation}
		where $ u_{t}^{\ast D} $ is the optimal solution to Problem D, and $ \nu_{t} $ is Lagrangian multiplier.
		
		As $ u_{t}^{\ast} $ and $ c_{t}^{\ast} $ are the optimal solutions to Problem B, if $ u_{t}^{\ast} = q_{t} $, according to (\ref{3.42}) and (\ref{3.43}), the following hold:
		\begin{equation}
		\begin{split}
		\label{3.55}
		\alpha_{t}^{\ast} \geq 0
		\end{split}
		\end{equation}
		\begin{equation}
		\begin{split}
		\label{3.56}
		\dfrac{\partial  F_{t}\left(q_{t},H_{t}^{\ast}\right) }{\partial u_{t}} =& \dfrac{-c^{\ast} + \alpha_{t}^{\ast}}{\delta_{t}} - \dfrac{\partial Y(H_{t+1}^{\ast})}{\partial u_{t}} 
		\\>& \dfrac{-c^{\ast} }{\delta_{t}} - \dfrac{\partial Y(H_{t+1}^{\ast})}{\partial u_{t}} 
		= -c^{\ast}_{t}
		\end{split}
		\end{equation}
		
		Given that $ F_{t}\left(u_{t},H_{t}^{\ast}\right) $ is convex, if $ u_{t}^{\ast D} \neq q_{t} $, then:
		\begin{equation}
		\begin{split}
		\label{3.57}
		\dfrac
		{\partial  F_{t}\left(u_{t}^{\ast D},H_{t}^{\ast}\right) }
		{\partial u_{t}} 
		\geq
		\dfrac
		{\partial  F_{t}\left(q_{t},,H_{t}^{\ast}\right) }
		{\partial u_{t}} 
		> -c^{\ast}_{t}
		\end{split}
		\end{equation}
		
		Combining equations (\ref{3.51}), (\ref{3.52}), and (\ref{3.57}), the following hold:
		\begin{equation}
		\begin{split}
		\label{3.58}
		\nu_{t}^{\ast} > 0
		\end{split}
		\end{equation}
		\begin{equation}
		\begin{split}
		\label{3.59}
		u_{t}^{\ast D} = q_{t}
		\end{split}
		\end{equation}
		
		So $ u_{t}^{\ast D} = u_{t}^{\ast} $. If $ u_{t}^{\ast} \geq q_{t} $, according to Theorem \ref{theorem1}, the following holds:
		\begin{equation}
		\begin{split}
		\label{3.60}
		\dfrac
		{\partial  F_{t}\left(u_{t}^{\ast},H_{t}^{\ast}\right) }
		{\partial u_{t}}  = -c^{\ast}_{t}
		\end{split}
		\end{equation}
		
		Then $ u_{t}^{\ast D} = u_{t}^{\ast} $ and $ \nu_{t}^{\ast} = 0 $ satisfy the KKT conditions of Problem D (\ref{3.51})-(\ref{3.54}). Given that $ F_{t}\left(u_{t},H_{t}\right) $ is convex with respect to $ u_{t} $, Problem D is also convex (as observed from (50)). Thus, $ u_{t}^{\ast} $ is the unique optimal solution to Problem D. When implementing $ u_{t} = u_{t}^{\ast} $ into Problem D (\ref{3.46}), it becomes:
		\begin{equation}
		\label{3.61}
		\begin{split}
		&\min\limits_{\bm{x}_{t}} \left[ f_{t}(\bm{x}_{t}) + c_{t}^{\ast}u_{t}^{\ast}\right] 
		\\ \mathrm{s.t.} \quad &\bm{g}_{t}(\bm{x}_{t},H_{t}^{\ast}) \leq \mathbf{0}
		\\ &d(\bm{x}_{t}^{\mathrm{S}},H_{t}^{\ast}) \leq u_{t}^{\ast}
		\\ &\bm{x}_{t}^{\mathrm{S}} \geq \mathbf{0}
		\end{split}
		\end{equation}
		which has the same solution set with Problem A', because the only difference between problems (\ref{3.19}) and (\ref{3.61}) is a constant in the objective function.
		
		Finally, as Problem D has the same solution set with both Problem A' and C, Lemma \ref{lemma2} is proved.	
	\end{proof}
	
	Combining Lemma \ref{lemma1} and \ref{lemma2}, Theorem \ref{theorem1} is straightforward on conditions that $F_{t}\left(u_{t},H_{t} \right)$ is differentiable and convex. When $F_{t}\left(u_{t},H_{t} \right)$ is only piece-wise differentiable and convex, it can be proved that for each subsets of $ u_{t} $ where $ F_{t}(u_{t},H_{t}) $ is  differentiable and convex, Theorem \ref{theorem1} holds, and so does the universal set of $ \bm{x}_{t} $.
	
	\subsection{Proof of Corollary \ref{corollary1}}
	\begin{proof}
		The Lagrangian function of Problem B is:
		\begin{equation}
		\begin{split}
		\label{3.62}
		L = &\sum\limits_{t \leq T} \delta_{t} F_{t}\left(u_{t},H_{t}\right) 
		+ c\left[ \sum\limits_{t \leq T}  u_{t} - U(H_{0}) \right] 
		\\+& \sum\limits_{t \leq T} \alpha_{t}(q_{t} - u_{t})
		\end{split}
		\end{equation}
		
		Then, two of the KKT conditions of Problem B are:
		\begin{equation}
		\begin{split}
		\label{3.63}
		\dfrac{\partial L}{\partial u_{t}} &= 
		\sum\limits_{t \leq \tau \leq T}\delta_{\tau} \dfrac{\partial F_{\tau}\left(u_{\tau}^{\ast},H_{\tau}^{\ast} \right)}{\partial u_{t}}
		+ c^{\ast} - \alpha_{t}^{\ast} = 0
		\end{split}
		\end{equation}
		\begin{equation}
		\begin{split}
		\label{3.64}
		\dfrac{\partial L}{\partial u_{t+1}} &= 
		\sum\limits_{t+1 \leq \tau \leq T}\delta_{\tau} \dfrac{\partial F_{\tau}\left(u_{\tau}^{\ast},H_{\tau}^{\ast} \right)}{\partial u_{t+1}}
		+ c^{\ast} - \alpha_{t+1}^{\ast} = 0
		\end{split}
		\end{equation}
		
		If equation (\ref{3.21}) is false, then:
		\begin{equation}
		\begin{split}
		\label{3.65}
		\alpha_{t}^{\ast} = \alpha_{t+1}^{\ast} = 0
		\end{split}
		\end{equation}
		
		For any $ \tau \geq t+1 $, the following holds:
		\begin{equation}
		\label{3.66}
		\begin{split}
		\dfrac{\partial H_{\tau} }{\partial u_{t}} = \dfrac{\partial H_{\tau} }{\partial u_{t+1}} = -\dfrac{H_{0} - \underline{H} }{ U(H_{0})} 
		\end{split}
		\end{equation}
		
		Subtracting equation (\ref{3.64}) from  equation (\ref{3.63}) and combining  equations (\ref{3.65}) and(\ref{3.66}), the following holds:
		\begin{equation}
		\begin{split}
		\label{3.67}
		\delta_{t} \dfrac{\partial F_{t}\left(u_{t}^{\ast},H_{t}^{\ast} \right)}{\partial u_{t}}
		+&
		\delta_{t+1} \dfrac{\partial F_{t+1}\left(u_{t+1}^{\ast},H_{t+1}^{\ast} \right)}{\partial u_{t}}
		\\-&
		\delta_{t+1} \dfrac{\partial F_{t+1}\left(u_{t+1}^{\ast},H_{t+1}^{\ast} \right)}{\partial u_{t+1}} = 0
		\end{split}
		\end{equation}
		
		Equation (\ref{3.67}) can be reformulated as: 
		\begin{equation}
		\begin{split}
		\label{3.68}
		c^{\ast}_{t} &= 
		- \dfrac{\partial F_{t}\left(u_{t}^{\ast},H_{t}^{\ast} \right)}{\partial u_{t}}
		\\&=\dfrac{\delta_{t+1}}{\delta_{t}}
		\left[ 
		- \dfrac{\partial F_{t+1}\left(u_{t+1}^{\ast},H_{t+1}^{\ast} \right)}{\partial u_{t}}
		+
		\dfrac{\partial F_{t+1}\left(u_{t+1}^{\ast},H_{t+1}^{\ast} \right)}{\partial u_{t+1}}
		\right] 
		\\&=\dfrac{\delta_{t+1}}{\delta_{t}}
		\left[ 
		c^{\ast}_{t+1}
		+
		\dfrac{\partial F_{t+1}\left(u_{t+1}^{\ast},H_{t+1}^{\ast} \right)}{\partial H_{t+1}}
		\dfrac{\partial H_{t+1}^{\ast}}{\partial u_{t+1}}
		\right] 
		\\&=\dfrac{\delta_{t+1}}{\delta_{t}}
		\left[ 
		c^{\ast}_{t+1}
		-
		\dfrac{\partial F_{t+1}\left(u_{t+1}^{\ast},H_{t+1}^{\ast} \right)}{\partial H_{t+1}}
		\dfrac{H_{0} - \underline{H} }{ U(H_{0})} 
		\right] 
		\end{split}
		\end{equation}
		which is exactly equation (\ref{3.22}).
	\end{proof}

	%
	%

	\ifCLASSOPTIONcaptionsoff
	\newpage
	\fi

	
	
	\bibliographystyle{IEEEtran}
	%
	\bibliography{Dissertation_Library}
	
	%
	
	
	
	
	
	
	

\end{document}